\documentclass[letterpaper, 11pt,  reqno]{amsart}

\usepackage{amsmath,amssymb,amscd,amsthm,amsxtra, esint}

\usepackage[margin=1.2in,marginparwidth=1.5cm, marginparsep=0.5cm]{geometry}

\usepackage[implicit=true]{hyperref}
\allowdisplaybreaks[2]

\usepackage{color}

\definecolor{gr}{rgb}   {0.,   0.69,   0.23 }
\definecolor{bl}{rgb}   {0.,   0.5,   1. }
\definecolor{mg}{rgb}   {0.85,  0.,    0.85}
\definecolor{yl}{rgb}   {0.8,  0.7,   0.}
\definecolor{or}{rgb}  {0.7,0.2,0.2}

\newtheorem{theorem}{Theorem} [section]

\newtheorem{lemma}[theorem]{Lemma}

\newtheorem{remark}[theorem]{Remark}

\newtheorem{definition}[theorem]{Definition}

\newcommand{\I}{\hspace{0.5mm}\text{I}\hspace{0.5mm}}

\newcommand{\noi}{\noindent}
\newcommand{\Z}{\mathbb{Z}}
\newcommand{\R}{\mathbb{R}}
\newcommand{\C}{\mathcal{C}}
\newcommand{\T}{\mathbb{T}}

\let\P= \undefined
\newcommand{\P}{\mathbf{P}}

\newcommand{\al}{\alpha}

\newcommand{\dl}{\delta}

\newcommand{\Dl}{\Delta}
\newcommand{\eps}{\varepsilon}

\newcommand{\ld}{\lambda}

\newcommand{\s}{\sigma}

\newcommand{\ft}{\widehat}

\newcommand{\wt}{\widetilde}

\newcommand{\dt}{\partial_t}
\newcommand{\dd}{\partial}

\newcommand{\les}{\lesssim}
\newcommand{\ges}{\gtrsim}

\newcommand{\M}{\mathcal{M}}

\newcommand{\N}{\mathbb{N}}

\newcommand{\too}{\longrightarrow}

\renewcommand{\I}{\mathcal{I}}

\newcommand{\TT}{\mathcal{T}}

\newcommand{\BT}{{\bf T}}

\usepackage{tikz}

\usetikzlibrary{shapes.misc}
\usetikzlibrary{shapes.symbols}
\usetikzlibrary{shapes.geometric}
\tikzset{
	dot/.style={circle,fill=black,draw=black,inner sep=0pt,minimum size=0.5mm},
	>=stealth,
	}
\tikzset{
	ddot/.style={circle,fill=white,draw=black,inner sep=0pt,minimum size=0.8mm},
	>=stealth,
	}

\tikzset{decision/.style={ 
        draw,
        diamond,
        aspect=1.5
    }}

\tikzset{dia2/.style
={diamond,fill=white,draw=black,inner sep=0pt,minimum size=1mm},
	>=stealth,
	}

\tikzset{dia/.style
={star,fill=black,draw=black,inner sep=0pt,minimum size=1mm},
	>=stealth,
	}

\makeatletter
\def\DeclareSymbol#1#2#3{\expandafter\gdef\csname MH@symb@#1\endcsname{\tikz[baseline=#2,scale=0.15]{#3}}}
\def\<#1>{\csname MH@symb@#1\endcsname}
\makeatother

\DeclareSymbol{X}{-2.4}{\node[dot] {};}
\DeclareSymbol{1}{0}{\draw[white] (-.4,0) -- (.4,0); \draw (0,0)  -- (0,1.2) node[dot] {};}
\DeclareSymbol{2}{0}{\draw (-0.5,1.2) node[dot] {} -- (0,0) -- (0.5,1.2) node[dot] {};}

\DeclareSymbol{1}{-2.7}
 {
  \draw (0,0) node[dot]{};}

\DeclareSymbol{1'}{-2.7}
 {\draw (0,0) node[ddot]{};}

\DeclareSymbol{1''}{-2.7}
 {\draw (0,0) node[dia]{};}

\DeclareSymbol{3}{0}
 {\draw (0,0) node[dot]{} -- (0,1.2) node[ddot] {}; 
 \draw (-1.2,0) node[dot] {} -- (0,1.2)node[ddot] {} -- (1.2,0) node[dot] {};}

\DeclareSymbol{3'}{0}
 {\draw (0,0) node[dia]{} -- (0,1.2) node[ddot] {}; 
 \draw (-1.2,0) node[dia] {} -- (0,1.2)node[ddot] {} -- (1.2,0) node[dia] {};}

\DeclareSymbol{31}{-3}{
\draw (0,-1)node[dot] {} -- (0,0) node[ddot] {}-- (0.9, -1) node[dot] {}; 
\draw (0,-1)node[dot] {} -- (0,0) node[ddot] {}-- (-0.9, -1) node[dot] {}; 
\draw (0,0)node[ddot] {} -- (1.3,1) node[ddot] {}-- (1.3, 0) node[dot] {}; 
\draw  (1.3,1) node[ddot] {}-- (2.6, 0) node[dot] {}; 
}

\DeclareSymbol{32}{-3}{
\draw (0,-1)node[dot] {} -- (0,0) node[ddot] {}-- (0.9, -1) node[dot] {}; 
\draw (0,-1)node[dot] {} -- (0,0) node[ddot] {}-- (-0.9, -1) node[dot] {}; 
\draw (0,0)node[ddot] {} -- (0,1) node[ddot] {}-- (0.9, 0) node[dot] {}; 
\draw (0,0)node[ddot] {} -- (0,1) node[ddot] {}-- (-0.9, 0) node[dot] {}; 
}

\DeclareSymbol{33}{-3}{
\draw (2.6,-1)node[dot] {} -- (2.6,0) node[ddot] {}-- (3.5, -1) node[dot] {}; 
\draw (2.6,-1)node[dot] {} -- (2.6,0) node[ddot] {}-- (1.7, -1) node[dot] {}; 
\draw (0,0)node[dot] {} -- (1.3,1) node[ddot] {}-- (1.3, 0) node[dot] {}; 
\draw  (1.3,1) node[ddot] {}-- (2.6, 0) node[ddot] {}; 
}

\DeclareSymbol{3''}{0}
 {\draw (0,0) node[dia]{} -- (0,1.2) node[ddot] {}; 
 \draw (-1.2,0) node[dot] {} -- (0,1.2)node[ddot] {} -- (1.2,0) node[dia] {};}

\DeclareSymbol{3'''}{0}
 {\draw (0,0) node[dot]{} -- (0,1.2) node[ddot] {}; 
 \draw (-1.2,0) node[dot] {} -- (0,1.2)node[ddot] {} -- (1.2,0) node[dia] {};}

\DeclareSymbol{31'}{-3}{
\draw (0,-1.2)node[dia] {} -- (0,0) node[ddot] {}-- (0.9, -1.2) node[dia] {}; 
\draw (0,-1.2)node[dia] {} -- (0,0) node[ddot] {}-- (-0.9, -1.2) node[dia] {}; 
\draw (0,0)node[ddot] {} -- (1.3,1.2) node[ddot] {}-- (1.3, 0) node[dia] {}; 
\draw  (1.3,1.2) node[ddot] {}-- (2.6, 0) node[dia] {}; 
}

\DeclareSymbol{32'}{-3}{
\draw (0,-1.2)node[dia] {} -- (0,0) node[ddot] {}-- (0.9, -1.2) node[dia] {}; 
\draw (0,-1.2)node[dia] {} -- (0,0) node[ddot] {}-- (-0.9, -1.2) node[dia] {}; 
\draw (0,0)node[ddot] {} -- (0,1.2) node[ddot] {}-- (0.9, 0) node[dia] {}; 
\draw (0,0)node[ddot] {} -- (0,1.2) node[ddot] {}-- (-0.9, 0) node[dot] {}; 
}

\DeclareSymbol{33'}{-3}{
\draw (2.6,-1.2)node[dia] {} -- (2.6,0) node[ddot] {}-- (3.5, -1.2) node[dia] {}; 
\draw (2.6,-1.2)node[dia] {} -- (2.6,0) node[ddot] {}-- (1.7, -1.2) node[dia] {}; 
\draw (0,0)node[dot] {} -- (1.3,1.2) node[ddot] {}-- (1.3, 0) node[dot] {}; 
\draw  (1.3,1.2) node[ddot] {}-- (2.6, 0) node[ddot] {}; 
}

\newtheorem*{ackno}{Acknowledgments}

\numberwithin{equation}{section}
\numberwithin{theorem}{section}

\makeatletter
\@namedef{subjclassname@2020}{%
  \textup{2020} Mathematics Subject Classification}
\makeatother

\begin{document}

\title[Norm inflation for the cubic nonlinear heat equation]
{Norm inflation for the cubic nonlinear heat equation above
the scaling critical regularity}

\author[I. Chevyrev, T. Oh, Y. Wang]
{Ilya Chevyrev, Tadahiro Oh, and Yuzhao Wang}

\address{
Ilya Chevyrev,  School of Mathematics\\
The University of Edinburgh\\
and The Maxwell Institute for the Mathematical Sciences\\
James Clerk Maxwell Building\\
The King's Buildings\\
Peter Guthrie Tait Road\\
Edinburgh\\ 
EH9 3FD\\
 United Kingdom}

\email{ichevyrev@gmail.com}

\address{
Tadahiro Oh, School of Mathematics\\
The University of Edinburgh\\
and The Maxwell Institute for the Mathematical Sciences\\
James Clerk Maxwell Building\\
The King's Buildings\\
Peter Guthrie Tait Road\\
Edinburgh\\ 
EH9 3FD\\
 United Kingdom,
 and
 School of Mathematics and Statistics, Beijing Institute of Technology,
Beijing 100081, China
}

\email{hiro.oh@ed.ac.uk}

\address{
Yuzhao Wang\\
School of Mathematics\\
Watson Building\\
University of Birmingham\\
Edgbaston\\
Birmingham\\
B15 2TT\\ United Kingdom}
\email{y.wang.14@bham.ac.uk}



\begin{abstract}
We 
consider  the ill-posedness issue for the cubic nonlinear heat equation 
and  prove norm inflation with infinite loss of regularity
 in the H\"older-Besov space $\C^s = B^{s}_{\infty, \infty}$
for $ s \le  -\frac 23$.
In particular, 
our result includes the subcritical  range
  $-1<  s \le -\frac 23$, which is above the scaling critical regularity $s = -1$
with respect to the H\"older-Besov scale.
In view of the well-posedness result in $\C^s$, $s > -\frac 23$, our ill-posedness
result is sharp.

\end{abstract}


\maketitle

\tableofcontents

\section{Introduction}

We consider 
the  cubic  nonlinear heat equation (NLH): 
\begin{align}
\begin{cases}
 \dt u -    \Dl u +  u^3  = 0\\
u|_{t = 0} = u_0, 
\end{cases}
\ (x, t) \in \M\times \R_+, 
\label{NLH1}
\end{align}

\noi
where $\M = \R^d$ or $\T^d$ with $\T^d = \R^d/(2\pi \Z)^d$.
The equation \eqref{NLH1} and  its variants (see, for example,  \eqref{NLH1a} and \eqref{NLH1b} below) appear
in various physical settings and have been studied
extensively
from both theoretical and applied points of view.
See \cite{W81, Giga, BC, R} for the classical well-posedness
results on \eqref{NLH1}.
Over the last two decades, there has been an increasing interest in 
studying the equation \eqref{NLH1} with 
(singular) stochastic forcing and\,/\,or random initial data; see for example \cite{DD, HRW, Hairer,HLR}.
Such study often requires the understanding of solutions to \eqref{NLH1}
in the low regularity setting, 
where the equation \eqref{NLH1} is believed to be ill-posed.
However, the aforementioned ill-posedness seems to be part of the folklore, 
and our main goal in this paper is to provide
a rigorous proof of this ill-posedness by exhibiting
a strong instability known as a {\it norm inflation}, 
thus completing the well-posedness study of~\eqref{NLH1}.

In the following, we study the equation \eqref{NLH1} with respect to the H\"older-Besov space
$\C^s(\M) = B^s_{\infty, \infty}(\M)$, where the latter space $B^s_{\infty, \infty}(\M)$ 
denotes the usual Besov space defined by the norm:
\begin{align}
\|u \|_{B^s_{\infty, \infty}} = \sup_{j \in \N \cup\{0\}} 2^{js} \| \P_j u \|_{L^\infty}.
\label{besov1}
\end{align}

\noi
Here, $\P_j$ denotes the non-homogeneous Littlewood-Paley projector
onto the frequencies $\{ \xi \in \ft \M: |\xi|\sim 2^j\}$ for $j \geq 1$
and  $\{ \xi \in \ft \M
: |\xi|\les  1\}$ for $j = 0$, where 
  $\ft \M$ denotes the Pontryagin dual of $\M$, 
i.e.~$\ft \M = \R^d$ if $\M = \R^d$ and $\ft \M = \Z^d$ if $\M  = \T^d$.
Recall the following  Schauder estimate:
\begin{align}
\|e^{t\Dl} f \|_{\C^{\s}(\M)}
\les t^{- \frac{\s - s}{4}}
\|e^{t\Dl} f \|_{\C^s (\M)}
\label{Sch1}
\end{align}

\noi
for any $\s \ge s$ and $0 < t \le 1$.\footnote{The upper bound on $t$ is needed on $\T^d$.}
When $\M = \R^d$, \eqref{Sch1} follows from 
Young's inequality and estimating 
the $W^{{\s-s}, 1}$-norm of the heat kernel  for $e^{t\Dl}$;
see \cite[Lemma 2.4]{BCD} (for the case $\s = s$).
When $\M = \T^d$, 
\eqref{Sch1}
follows from Young's inequality
and the Poisson summation formula 
to pass an estimate on (fractional derivatives of) the heat kernel on $\T^d$ to that in a weighted Lebesgue space on $\R^d$; see the proof of 
Lemma 2.3 in \cite{LiuOh}
for such a computation (but done for $e^{-t\sqrt{-\Dl}}$).
Then, 
using the Schauder estimate~\eqref{Sch1}, 
one can easily show that~\eqref{NLH1} is locally well-posed
in $\C^s(\M)$ for $s > -\frac 23$.
The well-posedness issue of \eqref{NLH1} for $s\le -\frac 23$, however, 
appears unresolved.

As is well known, symmetries provide important heuristics
in studying well-posedness issues of nonlinear PDEs.
We now recall the scaling symmetry for \eqref{NLH1};
 if $u(x, t)$ is a solution to~\eqref{NLH1}
on $\R^d$, then
$u^\ld(x, t) := \ld^{-1} u (\ld^{-1}x, \ld^{-2}t)$
is also a solution to~\eqref{NLH1} on $\R^d$ with the scaled initial data $u^\ld_0 := \ld^{-1} u (\ld^{-1}x)$.
A direct computation shows that
this scaling symmetry leaves the homogeneous Besov norm  $\|u_0^\ld\|_{\dot B^{-1}_{\infty, \infty}}$
invariant, thus inducing
the scaling critical regularity $s_\text{crit} := -1$
with respect to the H\"older-Besov scale.
 In studying an evolution equation with initial data
 in a Banach space $B^s$ of regularity $s$ (think of $B^s$ as $\C^s$, the Sobolev space $H^s$, etc.),
 it is commonly conjectured that the equation is 
 well-posed in $B^s$   for  $s > s_\text{crit}$,
 while it is ill-posed  for $s < s_\text{crit}$.
In fact, on function spaces of positive regularities, 
this conjecture has been confirmed affirmatively for many parabolic and dispersive PDEs.
However, working on function spaces of negative regularities, 
this heuristics provided by the scaling argument is known to break down 
in some instances; see \cite{IO, MT, Ok, CP, Kishimoto,  FO}
for ill-posedness results above critical regularities.
See also \cite{BP, CD, IO, Oh2, Oh, OOT} for ill-posedness at critical regularities.

In this paper, we establish  a norm inflation with infinite
loss of regularity for \eqref{NLH1} {\it above} the scaling critical regularity $s_\text{crit} = -1$.

\begin{theorem}\label{THM:1}
Given $d \in \N$, 
let  $\M = \R^d$ or $\T^d$.
Let $s \le -\frac 23$.
Then, 
given any $\eps > 0$, 
there exist a 
solution $u_\eps$ to \eqref{NLH1} on $\M$
and $t_\eps \in (0, \eps) $ such that 
\begin{align}
 \| u_\eps(0)  \|_{\C^s(\M)} < \eps \qquad \text{ and } 
\qquad \| \P_0 u_\eps(t_\eps)\|_{L^\infty(\M)} > \eps^{-1}.
\label{thm1}
\end{align}
\end{theorem}
 
H{\"o}lder--Besov spaces play a dominant role in parabolic singular SPDEs. The threshold $s > -\frac 23$ appears in the study of the parabolic $\Phi^4_3$ model \cite{CC18,Hairer} (see also \cite{HLR}), and our results prove
that  this threshold is optimal. 

We also remark that, since $\|u_\eps(t_\eps)\|_{\C^\s} \ge  \| \P_0 u_\eps(t_\eps)\|_{L^\infty}$, the second estimate in~\eqref{thm1}
implies $\|u_\eps(t_\eps)\|_{\C^\s} > \eps^{-1}$ for any $\s\in\R$.
When $\s = s$, the statement corresponds to the usual norm inflation
 introduced in \cite{CCT2b} in the context of 
 the nonlinear Schr\"odinger equations (NLS)
 and the nonlinear wave equations.
  In proving Theorem \ref{THM:1}, 
 we exploit  a robust high-to-low energy transfer,
 which allows us to prove \eqref{thm1}, thus
 leading to the so-called
 ``infinite loss of regularity'';
 see, for example, \cite{CK, FO}.\footnote{The work \cite{FO} does not state
 their result as infinite loss of regularity but it follows from the proof.}
 
 Note that 
  Theorem \ref{THM:1} in particular  implies discontinuity
  of the solution
 map   $\Phi: u_0\in  \C^s(\M) \mapsto u \in C([0, T]; \C^s(\M))$ 
 at the trivial function $u_0 \equiv 0$
 for $s \le -\frac 23$, thus establishing 
 ill-posedness of the cubic NLH \eqref{NLH1}.
Furthermore,  in view of the aforementioned local well-posedness in $\C^s(\M)$ for $s > -\frac 23$, 
 Theorem \ref{THM:1} is sharp.  
See also Remark \ref{rem:q}\,(iii).

Our proof of Theorem \ref{THM:1} is based on the Fourier analytic approach
\cite{IO, Kishimoto, Oh2}.
In~\cite{IO},   Iwabuchi and Ogawa
introduced a new method for
 proving ill-posedness of 
evolution equations, exploiting high-to-low energy transfer
in the second Picard iterate.
In \cite{Kishimoto}, Kishimoto systematically
applied this approach 
to prove ill-posedness of NLS.
See  the work by Bejenaru and Tao \cite{BT}
for a precursor of this approach.  
We also mention the work by Bourgain and Pavlov\'ic \cite{BP}
on the Navier-Stokes equations, where an analogous high-to-low energy transfer was exploited.
See also \cite{Yoneda, Wang}.
In \cite{IO, Kishimoto}, 
the (scaled) modulation space $M_{2, 1}$
and its algebra property played an important role.
In \cite{Oh2}, the second author
implemented a slight simplification of  this approach
by working on the Wiener algebra (instead of the modulation spaces).
We point out that the H\"older-Besov space $\C^s(\M)$ is 
not a Fourier lattice\footnote{A Fourier lattice is a space, 
where a norm depends only on the absolute value of the Fourier transform.}
and thus is not quite suitable for this Fourier analytic approach \cite{IO, Kishimoto, Oh2}.
For example, in \cite{IO, MT}, this method was applied 
to prove ill-posedness of the (fractional) heat equations but the results were restricted to the $L^2$-based Sobolev spaces $H^s$
and the $L^2$-based Besov spaces $B^s_{2, q}$, both of which are Fourier lattices.
As we see below, our choice of initial data turns out to be rather simple (see \eqref{data} and \eqref{Z1} below), 
enabling us to adapt this Fourier analytic approach to the $L^\infty$-based function spaces similar in spirit to \cite{BP}.

As in \cite{IO, Kishimoto, Oh2}, 
we first express a solution $u$ to \eqref{NLH1} with $u|_{t = 0} = u_0$ in the  power series expansion:
\begin{align*}
 u & 
 = \sum_{j = 0}^\infty \Xi_j (u_0),
 \end{align*}

\noi
where $\Xi_j(u_0)$ denotes 
homogeneous multilinear terms (in $u_0$) of 
degree $2j+1$. In \cite{IO, Kishimoto}, 
 $\Xi_j$'s were defined by a recursive relation.
In proving Theorem \ref{THM:1}, 
we instead follow the approach in \cite{Oh2}
and  define the $j$th term $\Xi_j$ directly
via 
the power series expansion indexed by trees.
See Section \ref{SEC:2}.
There are several advantages to this latter approach.
First, it  allows us to establish nonlinear estimates without an induction.
By  a slight modification, 
this approach allows us 
to treat the Allen-Cahn equation (see \eqref{NLH1a} below)
in a  straightforward manner.
Furthermore, arguing as in \cite{Oh2}, we can easily extend
Theorem~\ref{THM:1}
to a
norm inflation based at a general initial condition $u_0 \in \C^s(\M)$;
given any $u_0 \in \C^s(\M)$ and any $\eps > 0$,
there exist a solution $u_\eps$ to \eqref{NLH1}
and $t_\eps  \in (0, \eps) $ such that 
\begin{align}
 \| u_\eps(0) - u_0\|_{\C^s(\M)} < \eps \qquad \text{ and } \qquad
  \| \P_0 u_\eps(t_\eps)\|_{L^\infty(\M)} > \eps^{-1}, 
\label{general1}
 \end{align}

\noi
thus establishing discontinuity of the solution map
everywhere in $\C^s(\M)$.

\medskip

We conclude this introduction by several remarks.

\begin{remark}\rm
(i) 
Unlike the dispersive case treated in \cite{IO, Kishimoto, Oh2, CP}, 
we work in the $L^\infty$-based function spaces
and thus our argument is insensitive to 
dimensions $d \geq 1$.

\smallskip

\noi
(ii)
For the standard cubic nonlinearity, Choffrut and Pocovnicu \cite{CP}
proved a norm inflation 
for the fractional  nonlinear Schr\"odinger equation:
\[ i \dt u + (-\Dl)^\frac{\al}{2} u + |u|^2 u = 0\]

\noi
above the scaling critical regularity, provided that $\al > 2$.
Namely, a  dispersion stronger than the standard Laplacian was needed.
It is interesting to compare this with the parabolic case (Theorem \ref{THM:1}), 
where a norm inflation above the scaling critical regularity
hold true for the standard Laplacian.
See also \cite{Kishimoto} for a norm inflation above the scaling critical regularity
for the quadratic NLS with the less physical nonlinearity $|u|^2$.

\end{remark}

\begin{remark}\rm
(i)  While  we present the result for the defocusing case (i.e.~with the $+$-sign in~\eqref{NLH1}), 
 our argument is insensitive to the defocusing\,/\,focusing nature of the nonlinearity and thus it also applies to the focusing case:
\begin{align}
\dt u - \Dl u - u^3 = 0.
\label{eq:reverse_sign}
\end{align}
 
 \noi
Furthermore, our argument is robust so that it can be applied to treat 
 a general power-type nonlinearity $u^k$.
See \cite{Kishimoto,  FO} in the case of the nonlinear Schr\"odinger equations
 and the nonlinear wave equations with 
 a general power-type nonlinearity. 
 
 \smallskip
 
 \noi
 (ii) As mentioned above, a slight modification of the proof of Theorem 
 \ref{THM:1}
 yields a norm inflation with a general initial data as stated in \eqref{general1}.
 See \cite{Oh2} for such an argument.
\end{remark}

\begin{remark}\rm
Consider 
the 
Allen-Cahn equation:
\begin{align}
 \dt u -    \Dl u - u +  u^3  = 0.
 \label{NLH1a}
\end{align}

\noi
By the Schauder estimate, one can easily prove local well-posedness of \eqref{NLH1a}
in $\C^s(\M)$ for $s > -\frac 23$.
On the other hand, a slight modification of the proof of Theorem \ref{THM:1}
yields a norm inflation for \eqref{NLH1a} in $\C^s(\M)$ for $s \le - \frac 23$.
See Remark \ref{REM:AC}.

Next, consider 
the Cahn-Hilliard equation:
\begin{align}
 \dt u +   \Dl^2 u - \Dl  u^3  = 0.
 \label{NLH1b}
\end{align}

\noi
and its variants:
\begin{align}
 \dt u +   \Dl^2 u -   u^2 \Dl u   = 0
\qquad \text{or} \qquad  \dt u +   \Dl^2 u -   u (\dd u)(\dd  u)   = 0.
 \label{NLH1c}
\end{align}

\noi
These equations all scale the same way and the scaling critical regularity is $s = -1$
in the H\"older-Besov scale.
By  the Schauder estimate:
\begin{align}
\| e^{-t\Dl^2} u_0\|_{\C^\s} \les t^{- \frac{\s - s}{4}} \|u_0\|_{\C^s}
\label{S1}
\end{align}

\noi
for any $\s \ge s$ and $0 < t \le 1$,\footnote{As in \eqref{Sch1}, the upper bound on $t$ is needed on $\T^d$.}
we can prove local well-posedness of the Cahn-Hilliard equation~\eqref{NLH1b}
in $\C^s(\M)$ for the full subcritical range $s > -1$.
On the other hand, as for the equations in \eqref{NLH1c}, 
the Schauder estimate \eqref{S1} yields local well-posedness
in $\C^s(\M)$ only for a partial range $s > -\frac 23$.
In fact, for these variant  equations in \eqref{NLH1c}, 
our argument works and yields
 a norm inflation as in Theorem \ref{THM:1} holds  for $s \le -\frac 23$.
See Remark \ref{REM:AC}.

Due to the presence of the Laplacian on the cubic nonlinearity in \eqref{NLH1b}, 
our strategy for proving norm inflation via the high-to-low energy transfer 
(or using high-to-high energy transfer) does not work for 
the Cahn-Hilliard equation \eqref{NLH1b}.
It is of interest to note 
how the equations in \eqref{NLH1b} and \eqref{NLH1c} scale the same way
but that their well\,/\,ill-posedness results are quite different.
We point out that, as for proving a norm inflation, 
it seems easier to treat a nonlinearity that is not a total derivative, 
as in \eqref{NLH1c}.
See also   \cite{C, WZ}.

\end{remark}

\begin{remark}\label{rem:q}\rm
(i) It is worthwhile to note that the threshold regularity $s = -\frac 23 $ also appears in the study 
of long-time behavior of small data solutions.
In \cite{MYZ}, Miao, Yuan, and Zhang
proved small data global well-posedness of \eqref{NLH1} 
and \eqref{eq:reverse_sign} with small initial data in $\dot B^s_{p, q}(\R^d)$, 
where $s = \frac d p - 1$, $d < p < 3d$, and $1\le q \le \infty$.
Here, the regularity $s = \frac dp - 1$ corresponds to the scaling critical regularity 
for the $L^p$-based Besov spaces.
Note that the condition on $p$ implies that $- \frac 23 < s < 0$.
In \cite{BC2}, 
Brandolese and Cortez proved
existence of a (smooth) finite time blowup solution to \eqref{eq:reverse_sign}
with small initial data in $\dot B^{-\frac 23}_{3d, q}(\R^d)$, $3 < q \le \infty$,
thus providing a negative answer to Meyer's question \cite{Meyer}.

\smallskip

\noi
(ii) By a slight modification of the proof, 
we can extend the norm inflation to the $L^p$-based Besov spaces $B^s_{p, q}(\T^d)$
for any $1 \le p \le \infty$
in the case of the torus $\T^d$
for (a) $s < -\frac 23$ and $1 \le q \le \infty$
and (b) $s = -\frac 23$ and $3 < q \le \infty$.\footnote{In the endpoint case $s = -\frac 23$, 
we need to slightly modify the proof by replacing 
$\log K$ in \eqref{Z1} by $K^{\frac13 - \dl}$. The restriction on $q$ 
is needed to have  the  $B^s_{p, q}$-norm of the initial data tend to 0 as in \eqref{Z2}.}
On $\R^d$, we need to slightly modify our choice of initial data by localizing in space, 
but we expect the same result also holds on $\R^d$.

\smallskip

\noi
 (iii)
 Recall the following characterization of the Besov norm (\cite[Theorem 5.3]{LR}):
 \begin{align}
 \begin{split}
\|f\|_{B^{-\frac23}_{\infty,3}} 
& \sim \| e^\Dl f\|_{L^\infty}
+ \bigg(\int_0^1 t^{-1}\|t^{\frac1 3}e^{t\Delta} f\|_{L^\infty_x}^3dt \bigg)^{\frac 13} \\
& =
\| e^\Dl f\|_{L^\infty}
+ 
 \bigg(\int_0^1 \|e^{t\Delta} f\|_{L^\infty}^3dt \bigg)^{\frac 13}.
\end{split}
\label{AA1}
\end{align}

\noi
Using this characterization, 
it is straightforward to prove local well-posedness
of~\eqref{NLH1}  in $B^{-\frac23}_{\infty,3}(\T^d)$.
For example, one can first run a contraction argument
in $L^3([0, T]; L^\infty(\T^d))$  to construct a solution $u
\in L^3([0, T]; L^\infty(\T^d))$
 to \eqref{NLH1}, 
and then a posteriori show that $u$ belongs to $C([0, T]; B^{-\frac 23}_{\infty, 3}(\T^d))$.\footnote{Here,
as in \cite{MW1, MW2},  
we view the Besov space $B^{-\frac 23}_{\infty, 3}(\T^d)$
as the completion of $C^\infty(\T^d)$ under the norm~\eqref{AA1}.
Alternatively, i.e.~if we view the Besov space $B^{-\frac 23}_{\infty, 3}(\T^d)$ as the collection 
of elements in $\mathcal D'(\T^d)$ with  finite $B^{-\frac 23}_{\infty, 3}$-norms, 
then a solution $u$ to \eqref{NLH1} does not belong to 
$C([0, T]; B^{-\frac 23}_{\infty, 3}(\T^d))$ in general.
With this latter definition of the Besov space, even the linear solution $e^{t \Dl} u_0$
does not belong to 
$C([0, T]; B^{-\frac 23}_{\infty, 3}(\T^d))$ in general.
In this case, we can actually show that, given 
$u_0 \in B^{-\frac 23}_{\infty,3}(\T^d)$, 
 the nonlinear part 
$v(t) : = u(t) - e^{t\Delta} u_0$ belongs to 
 $C([0,T];L^\infty(\T^d))$ and 
 that the map $u_0 \mapsto v = v(u_0)$ is locally Lipschitz with respect to $u_0 \in B^{-\frac 23}_{\infty,3}(\T^d)$.
}
Therefore, 
our results (Theorem \ref{THM:1} and Remark \ref{rem:q}\,(ii) above) are  sharp across the entire Besov scale $B^s_{\infty,q}(\T^d)$ for $s\in\R$ and $1 \le q \le \infty$.

\end{remark}

\begin{remark}\label{REM:C}\rm
Using probabilistic methods, the first author~\cite{C} recently showed the same type of norm inflation as~\eqref{thm1} in the scaling subcritical space $\C^{-\frac 12}(\T^d)$ 
(but not in $B^{-\frac 12}_{\infty,q}(\T^d)$ for finite  $q<\infty$)
for a heat equation with dominant nonlinearity $u\times Du$ that is not a total derivative, 
where the scaling critical space is $\C^{-1}(\T^d)$.
\end{remark}

\section{Preliminary analysis}
\label{SEC:2}

In this section, we review a basic local well-posedness result for \eqref{NLH1}
and the power series expansion of a solution. While the following presentation closely follows
that in \cite{Oh2}, we include details for readers' convenience.

We say that $u$ is a solution to \eqref{NLH1}
with $u|_{t = 0} = u_0$
if $u$ satisfies the following Duhamel formulation:
\begin{align}
u(t) = e^{t\Dl} u_0 + \I[u](t), 
\label{Duhamel3}
\end{align}

\noi
where $\I$ denotes the 
 the Duhamel integral  operator
defined by 
\begin{align}
\I[u_1, u_2, u_3](t)
:= - 
\int_0^t e^{(t - t')\Dl}  \bigg(\prod_{j = 1}^3u_j (t') \bigg)dt'
\label{Duhamel1}
\end{align}

\noi
with a shorthand notation $\I[u] := \I[u, u, u]$
when all the three arguments $u_1, u_2$, and $u_3$ are identical.
Then, from the boundedness of the heat semigroup 
and the algebra property of 
 $\C^{s_0}(\M)$, $s_0 > 0$,  we have the following local well-posedness of \eqref{NLH1}
in $\C^{s_0}(\M)$, $s_0 > 0$.

\begin{lemma}\label{LEM:LWP}
Let $s_0 > 0$. Then, 
the cubic NLH \eqref{NLH1} is locally well-posed in $\C^{s_0}(\M)$.
More precisely, given $u_0 \in \C^{s_0}(\M)$, 
there exist $T \sim \| u_0\|_{\C^{s_0}}^{-2}>0 $ and 
a unique solution $ u \in C([0,  T]; \C^{s_0}(\M))$
satisfying \eqref{Duhamel3}.

\end{lemma}

%

Fix $u_0 \in \C^{s_0}(\M)$ for some $s_0 >0$,
and let $u$ be the solution to \eqref{NLH1} with $u|_{t=0} = u_0$
constructed in Lemma \ref{LEM:LWP}.
Following the approach in \cite{Oh2}, we express the solution $u$  in the power series expansion;
see  \cite{Sinai, Gub, CH2},  where the power series expansion was used to construct solutions
to parabolic and dispersive equations.
We first recall the following definition of (ternary) trees.

\begin{definition} \label{DEF:tree} \rm
(i) Given a partially ordered set $\TT$ with partial order $\leq$, 
we say that $b \in \TT$ 
with $b \leq a$ and $b \ne a$
is a child of $a \in \TT$,
if  $b\leq c \leq a$ implies
either $c = a$ or $c = b$.
If the latter condition holds, we also say that $a$ is the parent of $b$.

\smallskip 

\noi
(ii) 
A tree $\TT$ is a finite partially ordered set,  satisfying
the following properties:
\begin{itemize}
\item[(a)] Let $a_1, a_2, a_3, a_4 \in \TT$.
If $a_4 \leq a_2 \leq a_1$ and  
$a_4 \leq a_3 \leq a_1$, then we have $a_2\leq a_3$ or $a_3 \leq a_2$.

\item[(b)]
A node $a\in \TT$ is called terminal, if it has no child.
A non-terminal node $a\in \TT$ is a node 
with  exactly three children.

\item[(c)] There exists a maximal element $r \in \TT$ (called the root node) such that $a \leq r$ for all $a \in \TT$.

\item[(d)] $\TT$ consists of the disjoint union of $\TT^0$ and $\TT^\infty$,
where $\TT^0$ and $\TT^\infty$
denote  the collections of non-terminal nodes and terminal nodes, respectively.
\end{itemize}

\end{definition}

\noi
Note that the number $|\TT|$ of nodes in a tree $\TT$ is $3j+1$ for some $j \in \mathbb{N}\cup\{0\}$,
where $|\TT^0| = j$ and $|\TT^\infty| = 2j + 1$.
Let us denote  the collection of trees in the $j$th generation (i.e.~with $j$ parental nodes) by $\BT(j)$, i.e.
\begin{equation*}
\BT(j) := \big\{ \TT : \TT \text{ is a tree with } |\TT| = 3j+1 \big\}.
\end{equation*}

\noi
Then, there exists $C_0 > 0$ such that 
\begin{align}
\# \BT(j) \leq C_0^j
\label{tree0a}
\end{align}

\noi
for any $j \in \mathbb{N}\cup\{0\}$.
See \cite{Oh2} for the proof of \eqref{tree0a}.

We now define a map $\Psi = \Psi_{u_0}: 
\bigcup_{j = 0}^\infty \BT(j) \to \mathcal{D}'(\M\times [0, T])$
as follows.
Given a tree $\TT \in \BT(j)$,
$j \in \N \cup\{0\}$, 
we define $\Psi(\TT)$ by the following rules:
\begin{itemize}
\item[(i)] Replace a non-terminal node ``\,$\<1'>$\,'' 
by the Duhamel integral operator $\I$ defined in \eqref{Duhamel1}
with its three children as arguments $u_1, u_2$, and $u_3$.

\item[(ii)] Replace a terminal node ``\,$\<1>$\,'' 
by the linear solution $e^{t\Dl}u_0$. 
\end{itemize}

\noi
Note that, if  $\TT \in \BT(j)$, 
then $\Psi(\TT) $ is  $(2j+1)$-linear in $u_0$.
For example, we have
$\Psi(\,\<1>\,) = e^{t\Dl} u_0$, 
where  ``\,$\<1>$\,'' denotes the trivial tree, consisting only of the root node.
Similarly, we have 
\begin{align*}
\Psi( \<3>)  =
\I[e^{t\Dl} u_0]
\quad \text{and}\quad
\Psi\big( \<31>\big)  =
\I[\I[e^{t\Dl} u_0], e^{t\Dl} u_0, e^{t\Dl} u_0].
\end{align*}

\noi
Lastly, we define $\Xi_j$ by 
\begin{align}
\Xi_j (u_0)
: =  \sum_{\TT \in \BT(j)} \Psi(\TT).
\label{tree1}
 \end{align}

Denoting the solution $u$ by  a star-shaped terminal node ``$\,\<1''>\,$'',\footnote{More
precisely, we set $\Psi(\<1''>) = u$.}
we can express the Duhamel formulation \eqref{Duhamel3} as 
\begin{align}
 \<1''>\, = \, \<1>\, + \, \<3'>\,.
\label{R1}
 \end{align}

\noi
By recursively applying \eqref{R1}
and eliminating the occurrence of ``$\,\<1''>\,$''
from younger trees, we have
\begin{align}
\begin{split}
 \<1''>\, 
 & = \, \<1>\, + \, \<3''>\, + \,\<31'>\, \\
 & = \, \<1>\, + \, \<3'''> \,+\, \<31'> \,+\,   \<32'> \\
 & = \, \<1>\, + \, \<3> \,+\,\<31'>\,+ \,\<32'>\,+\,   \<33'> 
 \\
 &  = \cdots = \, \<1>\, + \, \<3> \,+\,\<31>\,+ \,\<32>\,+\,   \<33> 
 + \cdots.
 \end{split}
 \label{R2}
 \end{align}

\noi
By extending the definition of $\Psi$
to formal sums of trees via linearity
and applying $\Psi$ to~\eqref{R2}, 
we obtain  the following (formal) power series
expansion of the solution $u$ to \eqref{NLH1} with $u|_{t = 0} = u_0$:
\begin{align}
\begin{split}
 u & 
 = \sum_{j = 0}^\infty \Xi_j (u_0)
 = \sum_{j = 0}^\infty \sum_{\TT \in \BT(j)} \Psi(\TT)\\
& = 
\Psi(\, \<1>\,) 
+ \Psi( \<3>) +
\Psi\big( \<31>\big) 
+
\Psi\big( \<32>\big) 
+ \Psi\big( \<33>\big) 
+\cdots.
\end{split}
\label{R3}
 \end{align}

\noi
Then, from 
\eqref{R3} with 
\eqref{Duhamel1}, 
\eqref{tree0a}, \eqref{tree1}, 
the boundedness of the heat semigroup, 
and the algebra property of 
 $\C^{s_0}(\M)$, $s_0 > 0$,
we have the following convergence result.

\begin{lemma}\label{LEM:LWP2}
Let $s_0 > 0$. Then, given $u_0 \in \C^{s_0}(\M)$,
the power series expansion in \eqref{R3} converges
absolutely  in $C([0,  T]; \C^{s_0}(\M))$, 
provided that  $T \le c_0 \| u_0\|_{\C^{s_0}}^{-2} $ for some small $c_0 > 0$.

\end{lemma}

It is easy to check that $u$ defined by the power series in \eqref{R3} 
indeed satisfies the Duhamel formulation \eqref{Duhamel3}.
Note that the time $T\sim \| u_0\|_{\C^{s_0}}^{-2}>0 $  for the local well-posedness in Lemma~\ref{LEM:LWP}
and for the convergence of the power series in Lemma~\ref{LEM:LWP2}
can be chosen to be the same.
Thanks to the unconditional uniqueness of the solution in 
the class  $\C([0,  T]; \C^{s_0}(\M))$, 
we conclude that 
$u$ defined by the power series in \eqref{R3} coincides with 
the solution constructed in Lemma \ref{LEM:LWP}.

\section{Non-endpoint case: $s < -\frac 23$}
\label{SEC:proof}

In this section, we present the proof of Theorem \ref{THM:1} for $s < -\frac 23$.
We will treat the endpoint case $s= -\frac 23$ in the next section.
Given $s < -\frac 23$,  the $\C^s$-norm is controlled by 
the $\C^{-\frac 23}$-norm and thus 
Theorem \ref{THM:1} for the non-endpoint case $s < -\frac 23$
follows from that for the endpoint case 
$s= -\frac 23$ whose proof is presented in the next section.
For readers' convenience, however, 
we decided to include
 the proof of Theorem \ref{THM:1} for the non-endpoint case, 
 since the relevant argument  is much simpler than that for the endpoint case
and thus is easier to follow.
As in the previous work \cite{IO, Kishimoto, Oh2}, the main strategy is to choose a suitable initial condition $u_0 = u_0(N)$, 
depending on a large parameter $N \gg1$ and vanishing as $N \to \infty$, 
and show that the Picard second iterate $\Xi_1(u_0)$ diverges as $N \to \infty$, 
while keeping the other terms under control.
In order to achieve the growth of the Picard second iterate, 
we exploit the high$\times$high$\times$high-to-low energy transfer mechanism 
of the cubic nonlinearity.


Note that in proving Theorem \ref{THM:1}, the $\C^s$-norm appears
only 
in the first estimate of  \eqref{thm1}, where we need to establish an upper bound.
Thus, in view of  $\| u \|_{\C^\s} \le \| u \|_{\C^s}$ for $\s \le s$, 
it suffices to only consider the subcritical case $-1< s < -\frac 23$.
Fix a large integer $N \gg 1$.
Then, we set the initial condition $u_0 = u_0(N)$ as
\begin{align}
    \label{data}
    u_0 (x)= N^{-s-\eps} \big( \cos (n_1 \cdot x) + \cos (2 n_1 \cdot x) \big)
\end{align}

\noi 
for some small $\eps > 0$, 
where $n_1 = (N, 0,\ldots,0)\in\Z^d$.
Then, from \eqref{besov1}, we have 
\begin{align*}
    \| u_0\|_{\C^s}  \sim N^{-\eps}  \too 0,
\end{align*}

\noi 
as $N \to \infty$.  This shows the first claim in \eqref{thm1}.

From \eqref{data}, we can write the linear solution as 
\begin{equation}\label{eq:linear_sol}
\begin{split}
e^{t\Delta} u_0 (x) & = N^{-s-\eps} \big( e^{-t|n_1|^2} \cos (n_1 \cdot x) + e^{-t 4 |n_1|^2} \cos (2 n_1 \cdot x) \big)\\
& =  N^{-s -\eps} \big( e^{-tN^2} \cos (n_1 \cdot x) + e^{-4tN^2} \cos (2n_1 \cdot x) \big),
\end{split}
\end{equation}

\noi
and thus we have 
\begin{align}
\P_0 e^{t\Delta} u_0 = 0.
\label{data2}
\end{align}

By a direct computation, we have
\begin{align*}
 (e^{t\Delta} u_0 (x) )^3 
 & = N^{-3s-3\eps} \big(  e^{-3tN^2} \cos^3 (n_1 \cdot x) + 3 e^{-6tN^2} \cos^2 (n_1 \cdot x)  \cos ( 2 n_1 \cdot x) \\
 & \quad + 3 e^{-9tN^2} \cos (n_1 \cdot x)  \cos^2 (2 n_1 \cdot x) + e^{-12tN^2} \cos^3 (2n_1 \cdot x)  \big). 
\end{align*}

\noi
Using the trigonometric identities:
\begin{align}
\cos^3 x &= \frac14 (\cos 3x + 3\cos x), \label{H1aa}\\
\cos(x_1) \cos (x_2) \cos (x_3)&  = \frac 14 \sum_{\eps_2, \eps_3 \in \{1, -1\}} \cos(x_1 + \eps_2  x_2 + \eps_3 x_3) ,\notag
\end{align}

\noi
we then have
\begin{align}
\begin{split}
 (e^{t\Delta} u_0 (x) )^3 & = \frac14 N^{-3s-3\eps} \Big[  e^{-3tN^2} (\cos (3 n_1 \cdot x) + 3\cos (n_1 \cdot x) \big) \\
 &\hphantom{XX}  + 3 e^{-6tN^2} \big(1 + \cos (4n_1 \cdot x) + 2 \cos ( 2 n_1\cdot x) \big) \\
 & \hphantom{XX} + 3 e^{-9tN^2} \big(\cos (5n_1 \cdot x) + 2\cos(n_1 \cdot x) + \cos (3n_1 \cdot x)  \big) \\
 & \hphantom{XX} + e^{-12tN^2} \big(\cos (6n_1 \cdot x) + 3 \cos (2n_1\cdot x) \big)   \Big] .
\end{split}
\label{Y1}
\end{align}

\noi
Hence, from \eqref{tree1} and \eqref{Y1}, we obtain
\begin{align*}
\Xi_1(u_0)(t) & = \int_0^t  e^{(t-t')\Delta} (e^{t'\Delta} u_0  )^3 dt' \\
& = \frac 18 N^{-2 -3s-3\eps} (1-e^{-6tN^2}) 
+ \frac14 N^{-3s-3\eps} \Big[   e^{-9tN^2} \frac{e^{6t N^2} -1}{6 N^2} \cos \big(3 n_1 \cdot x) \\
&  \hphantom{X}+ 3 e^{-tN^2} \frac{1- e^{-2tN^2}}{2N^2} \cos (n_1 \cdot x) \big) 
  + 3 e^{-16tN^2} \frac{e^{10tN^2} -1}{10N^2} \cos (4n_1 \cdot x)
\\
 & \hphantom{X}  + 3 e^{-4tN^2} \frac{1- e^{-2tN^2}}{N^2}  \cos ( 2 n_1\cdot x)
 + 3 e^{-25tN^2} \frac{e^{16tN^2} -1}{16N^2}\cos (5n_1 \cdot x)  \\
 & \hphantom{X}  + 3 e^{-tN^2} \frac{1-e^{-8tN^2}}{4N^2}  \cos(n_1 \cdot x) + 3 t e^{-9 tN^2} \cos (3n_1 \cdot x) \\
 & \hphantom{X}  + e^{-36tN^2} \frac{e^{24tN^2} -1}{24N^2} \cos (6n_1 \cdot x) + 3  e^{-4tN^2} \frac{1-e^{-8tN^2}}{8N^2}  \cos (2n_1\cdot x) \big)   \Big] .
\end{align*}

\noi
In particular, for $N \gg1$, we have
\begin{align*}
\P_0 \Xi_1(u_0)(t) 
& = \frac 18 N^{-2 -3s-3\eps} (1-e^{-6tN^2}) , 
\end{align*}

\noi
where $\P_0$ denotes the Littlewood-Paley frequency projector
onto the frequencies 
$\{ \xi \in \ft \M : |\xi|\les  1\}$.
Therefore, by  choosing $t = N^{-2+ 2\eps}$ such that $e^{-6tN^2}\ll 1$ for $N \gg 1$, 
we have 
\begin{align}
\| \P_0 \Xi_1(u_0)(t)\|_{L^\infty}
\sim  N^{-2 -3s-3\eps} \too \infty, 
\label{Y3}
\end{align}

\noi
as $N \to \infty$, provided that $s < -\frac 23 - \eps$.
Given $ s < -\frac 23$, 
this last condition can be guaranteed by choosing $\eps > 0$ sufficiently small.

Next, we estimate the contribution from the higher order terms in \eqref{R3}.
From
\eqref{tree1},  \eqref{tree0a}, 
the boundedness of the heat semigroup on $L^\infty(\M)$, 
and the algebra property of $L^\infty(\M)$ together with \eqref{data}
and $t = N^{-2+2 \eps}$,  
we have
\begin{align}
\begin{split}
\bigg\|\sum_{j = 2}^\infty  \P_0\Xi_j(u_0) (t) \bigg\|_{L^\infty}
& \les \sum_{j= 2}^\infty\| \Xi_j (u_0) (t) \|_{L^{\infty} }\\
& \le \sum_{j = 2}^\infty C_0^j t^j \| u_0 \|_{L^\infty}^{2j+1} \\
& \les \sum_{j =  2}^\infty C_0^j N^{(-2  - 2s ) j - s - \eps} \\
& \les N^{ - 4 -5s  - \eps}
\end{split}
\label{higher}
\end{align}

\noi
for $N \gg 1$, 
provided that $-2  - 2s  < 0$, namely 
$s > -1$.
Here, we used  that, by Lemmas \ref{LEM:LWP} and \ref{LEM:LWP2} with \eqref{data}, it suffices to have
$t = N^{-2+2\eps} \les \|u_0\|_{\C^{s_0}}^{-2} \sim N^{2s + 2\eps - 2s_0}$ for some $s_0>0$ to guarantee local well-posedness
and convergence of the power series expansion \eqref{R3}.\footnote{Here, we need 
$s_0>0$ such that $\C^{s_0}(\M)$ is an algebra, which is crucial for Lemmas \ref{LEM:LWP} and \ref{LEM:LWP2}.}
This condition reduces to 
 $s \geq  - 1 + s_0$, which can be guaranteed
 for given $s > -1$ by choosing $s_0 > 0$ sufficiently small.
We point out that 
the right-hand side of \eqref{higher} is controlled by \eqref{Y3}, provided that 
$s > -1 +  \eps$.
Given $s > -1$, 
this last condition can be guaranteed by choosing $\eps > 0$ sufficiently small.

Finally, putting 
 \eqref{R3}, 
\eqref{data2},
\eqref{Y3},
and \eqref{higher} together, we have
\begin{align}
    \label{final}
    \begin{split}
\| \P_0 u (t) \|_{L^\infty} 
& \ge \|\P_0\Xi_1(u_0)(t)\|_{L^\infty}
  - \| \P_0 e^{t\Delta} u_0\|_{L^\infty}-  \bigg\| \sum_{j = 2}^\infty \P_0\Xi_j (t) \bigg\|_{L^\infty}\\
& \ges N^{-2 -3s-3\eps} \too \infty, 
    \end{split}
\end{align}

\noi
as $N\to \infty$, provided that
\begin{align}
- 1+\eps < s< - \frac23-\eps
\label{Y4}
\end{align}

\noi
and $t = N^{-2 + 2\eps}$.
This proves Theorem \ref{THM:1} for $ s < -\frac 23$.

\begin{remark}\label{REM:AC}\rm
In this remark, we indicate how to modify the argument to treat the 
Allen-Cahn equation \eqref{NLH1a} and the 
Cahn-Hilliard-type equation \eqref{NLH1c}.
We only discuss the case $-1 < s < -\frac 23$.

We first consider the Allen-Cahn equation \eqref{NLH1a}.
Note that the difference 
between the heat semigroup and  the linear semigroup
$e^{t(1+\Dl)}$ for the Allen-Cahn equation appears only in the low frequencies.
Moreover, in proving a norm inflation, we study the dynamics only for very short times
and thus their difference becomes negligible.  Hence,  the argument presented above
can be applied directly to \eqref{NLH1a}.
In the following, however, we present an alternative argument, 
where we view $u - u^3$ as a combined power-type nonlinearity
and modify the structure of trees to allow nodes with different numbers of children.
We decided to include this argument since it can be adapted to 
study the nonlinear heat equation with a combined power-type nonlinearity
(such as the cubic-quintic nonlinearity).
See also \cite{WZ}, where ternary-quinary trees
were used to prove a norm inflation for the derivative nonlinear Schr\"odinger equation.


By writing \eqref{NLH1a} in the Duhamel formulation, we have
\begin{align*}
u(t) = e^{t \Dl} + \I^1[u] + \I[u], 
\end{align*}

\noi
where $\I$ is as in \eqref{Duhamel1}
and 
\begin{align}
\I^1[u] := - \I[u, 1, 1].
\label{I1}
\end{align}
Thus, we need to modify Definition \ref{DEF:tree} to  take 
this extra term $\I^1[u]$ into account.
More precisely, we need to replace the condition (b) by 

\smallskip
\begin{itemize}
\item[(b')]
A node $a\in \TT$ is called terminal, if it has no child.
A non-terminal node $a\in \TT$ is a node 
with  exactly one child or three children. 
\end{itemize}

\smallskip

Given a tree $ \TT $ with this new definition, 
let $|\TT|_1 $ (and $|\TT|_3$, respectively)  denotes the numbers of nodes 
which have one child (and three children, respectively).
Then, for $j_1, j_3 \in \N \cup \{0\}$, 
we set 
\begin{align*}
 \BT^{1,3}(j_1, j_2) & := \big\{ \TT :  \TT \text{ is a tree with } (|\TT|_1, |\TT|_3)
 =(j_1,j_3) \}.
\end{align*}

\noi
Note that for $\TT \in \BT^{1, 3}(j_1, j_3)$, 
we have
$|\TT| = j_1 + 3j_3+1$.
Furthermore, arguing as in the proof of Lemma 2.3 in \cite{Oh2}, 
the following bound on the size of $ \BT^{1, 3}(j_1, j_3)$ holds; 
there exists $C_0 > 0$ such that 
\begin{align*}
\sum_{j = j_1 + j_3}\#  \BT^{1, 3}(j_1, j_3) \leq C_0^{j}
\end{align*}

\noi
for any $j\in \mathbb{N}\cup\{0\}$.

We define
a map
 $\Psi = \Psi_{u_0}: 
\bigcup_{j_1, j_3 = 0}^\infty \BT^{1, 3}(j_1, j_3) \to \mathcal{D}'(\M\times [0, T])$
by replacing the definition (i) by 
\begin{itemize}
\item[(i')] Replace a non-terminal node ``\,$\<1'>$\,'' with three children
by the Duhamel integral operator $\I$ defined in \eqref{Duhamel1}
with its three children as arguments $u_1, u_2$, and $u_3$, 
and replace a non-terminal node ``\,$\<1'>$\,'' with one child
by $\I^1$ defined in \eqref{I1} with its child as an argument $u$.
\end{itemize}

\noi
We then set 
\begin{align*}
\Xi_j (u_0)
: =  \sum_{\substack{\TT \in \BT^{1, 3}(j_1, j_3)\\j = j_1 +j_3}} \Psi(\TT).
 \end{align*}

\noi
Then, the power series expansion for a (smooth) solution $u$ to \eqref{NLH1a}
is given by 
\begin{align*}
 u & 
 = \sum_{j = 0}^\infty \Xi_j (u_0)
 = \sum_{j = 0}^\infty  \sum_{\substack{\TT \in \BT^{1, 3}(j_1, j_3)\\j = j_1 +j_3}}\Psi(\TT).
 \end{align*}

\noi
We point out that Lemmas \ref{LEM:LWP} and \ref{LEM:LWP2} also hold for \eqref{NLH1a}.
Then, with the bound \eqref{tree0a}, 
we can simply repeat the argument presented above
to conclude a norm inflation for \eqref{NLH1a}.

Next, we consider the Cahn-Hilliard-type equation \eqref{NLH1c}.
Let us consider the first equation.  A similar argument applies to the second
equation in \eqref{NLH1c}.
In this case, we need to use the following Duhamel integral operator:
\begin{align}
\wt \I[u_1, u_2, u_3](t)
:= 
\int_0^t e^{-(t - t')\Dl^2}  u_1(t') u_2(t') \Dl u_3 (t')dt'.
\label{H0}
\end{align}

With the same choice \eqref{data} for the initial data, 
by repeating the computation above, 
we have
\begin{align*}
\P_0 \wt \Xi_1(u_0)(t) 
& = c_0 N^{-2 -3s-3\eps} (1-e^{-6tN^4}), 
\end{align*}

\noi
where $\wt \Xi_j$ is defined as in \eqref{tree1} but with the Duhamel 
integral operator $\wt \I$ in defined \eqref{H0}.
Hence, by choosing $t =  N^{-4+2\eps}$, 
the growth \eqref{Y3} still holds for $s < -\frac 23 - \eps$.

On the other hand, 
we claim that 
\begin{align}
\|\wt \Xi_j(u_0) (t) \|_{L^\infty}
\le C_0^j t^\frac j2 \| u_0\|^{2j+1}_{L^\infty}.
\label{H1}
\end{align}

\noi
For now, assume \eqref{H1} and we bound the higher order terms.
With the choice of $t$ as above,  \eqref{H1}, and \eqref{data}, we  have
\begin{align*}
\bigg\|\sum_{j = 2}^\infty  \P_0\wt \Xi_j(u_0) (t) \bigg\|_{L^\infty}
& \le \sum_{j = 2}^\infty C_0^j t^\frac{j}{2} \| u_0 \|_{L^\infty}^{2j+1} \\
& \les \sum_{j =  2}^\infty C_0^j N^{(-4  - 2s ) j - s - \eps} \\
& \les N^{-4 -5s   - \eps}
\end{align*}

\noi
just as in \eqref{higher}, 
provided that $-4  - 2s  < 0$, namely 
$s > -2$.
Hence, a norm inflation holds for 
the first equation in \eqref{NLH1c}.

It now remains to prove \eqref{H1}.
With \eqref{tree1} (for $\wt \Xi_j(u_0)$) and \eqref{tree0a}, 
it suffices to prove
\begin{align}
\|\Psi(\TT)(t)\|_{L^\infty} \les 
 t^\frac j2 \| u_0\|^{2j+1}_{L^\infty}
\label{H2}
\end{align}

\noi
for each $\TT \in \BT(j)$. We prove \eqref{H2} by induction.
When $j = 1$, 
the bound \eqref{H2} follows from 
the Schauder estimate \eqref{S1} (but for $L^\infty(\M)$): 
\begin{align}
\|\Dl e^{-t\Dl ^2} f \|_{L^\infty} \les t^{-\frac12} \| f\|_{L^\infty}.
\label{H3}
\end{align}

\noi
Now, fix $j_0 \ge 2$ and assume that \eqref{H2} holds for any  $j \le j_0 - 1$.
Given  $\TT \in \BT(j_0)$, 
suppose that there exists a terminal node  $a \in \TT^\infty$ 
whose parent $b$  is not the third child of its parent.\footnote{Here, the third child 
corresponds to $u_3$ in \eqref{H0} at the level of the Duhamel integral operator.}
Let $\TT_a$ denote the sub-tree in $\TT$ of one generation, containing $a \in \TT^\infty$
as one of its terminal nodes and $b$ as the parental node.
Then, let $\TT'$ be the tree obtained from $\TT$ by replacing the sub-tree $\TT_a$
with a node $b$
(and thus $b$ becomes a terminal node).
Note that  $\TT$ can be constructed by replacing
the terminal node $b \in (\TT')^\infty$ by the tree $\TT_a$.
By applying the inductive hypothesis
twice, we have
\begin{align*}
\|\Psi(\TT)(t)\|_{L^\infty}
&  \les 
 t^\frac {j_0 - 1} 2 \| u_0\|^{2j_0 - 2}_{L^\infty}
 \| \Psi(\TT_a) \|_{L^\infty}\\
&  \les 
 t^\frac {j_0 } 2 \| u_0\|^{2j_0 +1}_{L^\infty}, 
\end{align*}

\noi
yielding \eqref{H2}  for $j = j_0$ in this case.
Hence, it remains to consider the case when  there is no terminal node
whose parent is not the third child of its parent.
In other words, every parent in this (unbalanced)  tree is the third child of its parent.
In this case, by  the Schauder estimate \eqref{H3}, we have

\begin{itemize}
\item[(i)] for the last generation:\footnote{Here, we used the fact that 
$\int_0^t (t - t')^{-\frac{1}{2}}(t')^{-\frac{1}{2}}dt' = B(\frac 12 , \frac 12)$, 
where $B(a, b)$ is the beta function.}
\begin{align*}
\| \Dl \wt \I[e^{-t \Dl^2} u_0 ]\|_{L^\infty} \les \| u_0\|_{L^\infty}^3, 
\end{align*}

\item[(ii)]
for the intermediate generations:
\begin{align*}
\| \Dl \wt \I'[e^{-t \Dl^2} u_0 , e^{-t \Dl^2}u_0, u_3 ]\|_{L^\infty} \les t^\frac 12 \| u_0\|_{L^\infty}^2
\| u_3 \|_{L^\infty}, 
\end{align*}

\item[(iii)]
for the first  generation:
\begin{align*}
\| \wt \I'[e^{-t \Dl^2}u_0, e^{-t \Dl^2}u_0, u_3 ]\|_{L^\infty} \les t\| u_0\|_{L^\infty}^2
\| u_3 \|_{L^\infty}, 
\end{align*}

\end{itemize}

\noi
where $\wt \I'$ is defined by 
\begin{align*}
\wt \I'[u_1, u_2, u_3](t)
:= 
\int_0^t e^{-(t - t')\Dl^2}  \bigg(\prod_{j = 1}^3 u_j (t')\bigg)dt'.
\end{align*}

\noi
Putting (i), (ii), and (iii) together, 
we directly obtain \eqref{H2} in this case (without using the inductive hypothesis).
Therefore, by induction, \eqref{H2} holds for any $j \in \N$.

\end{remark}

\section{Endpoint case: $s = -\frac 23$}
\label{SEC:proof2}

In this section, we present the proof of Theorem \ref{THM:1}
at the endpoint regularity $s = -\frac 23$.
When $s = -\frac 23$, the construction in the previous section (with $\eps = 0$)
only gives the $O(1)$ growth of the Picard second iterate.
In order to overcome this difficulty, we consider a lacunary sequence
of initial data of type \eqref{data}, where each component yields  $O(1)$ contribution
to the Picard second iterate at the frequency 0.

Fix  large integers $N, K\gg 1$.
Then, we set the initial condition $u_0 = u_0(N, K)$ as
\begin{align}
    u_0 (x)= \frac{1}{\log K}\sum_{k = 1}^K (a_k N)^{\frac 23} \big( \cos (a_k n_1 \cdot x) + \cos (2a_k n_1 \cdot x) \big),
\label{Z1}
\end{align}

\noi 
where $n_1 = (N, 0,\ldots,0)\in\Z^d$ as before and 
\begin{align}
a_k = 2^{2^k}.
\label{Z1a}
\end{align}

\noi
Then, from \eqref{besov1}, we have 
\begin{align}
    \| u_0\|_{\C^{-\frac 23}}  \sim \frac{1}{\log K}  \too 0,
\label{Z2}
\end{align}

\noi 
as $K \to \infty$.  This shows the first claim in \eqref{thm1}.

Proceeding as in~\eqref{eq:linear_sol} and~\eqref{data2} with \eqref{Z1}, we have 
\begin{align}
\P_0e^{t\Delta} u_0 (t) =0. 
\label{Z3}
\end{align}

\noi
By repeating the computation in the previous section and making use
of the lacunarity of the summands in \eqref{Z1}, we have
\begin{align*}
 \P_0 (e^{t\Delta} u_0  )^3 
 & = \frac{3}{4(\log K)^3} \sum_{k = 1}^ K  (a_k N)^{2} e^{-6t(a_k N)^2}.
\end{align*}

\noi
Hence, by  choosing $t = N^{-2+\dl}$ for some small $\dl > 0$, we have
\begin{align}
\begin{split}
\| \P_0 \Xi_1(u_0)(t)\|_{L^\infty} &= \bigg\|\int_0^t  e^{(t-t')\Delta} \P_0 (e^{t'\Delta} u_0  )^3 dt'\bigg\|_{L^\infty} \\
&= \frac {1}{8(\log K)^3} \sum_{k = 1}^K
(1-e^{-6t(a_k N)^2}) \\
& \ges \frac{K}{(\log K)^3} \too \infty, 
\end{split}
\label{Z5}
\end{align}

\noi
as $K \to \infty$.

It remains to estimate the contribution from the higher order terms in \eqref{R3}.
Proceeding as in the previous section  with \eqref{Z1}, \eqref{Z1a}, 
and $t = N^{-2+\dl}$ and choosing $K = c_0 \log \log N$ for some small $c_0 >0$, we have
\begin{align}
\begin{split}
\bigg\|\sum_{j = 2}^\infty \P_0 \Xi_j(u_0) (t) \bigg\|_{L^\infty}
& \les \sum_{j= 2}^\infty\| \Xi_j (u_0) (t) \|_{L^{\infty} } \le \sum_{j = 2}^\infty C_0^j t^j \| u_0 \|_{L^\infty}^{2j+1} \\
& \les \sum_{j =  2}^\infty \frac{C_0^jK^{2j+1}}{(\log K)^{2j+1}} 
2^{(2j+1)2^K}N^{(-\frac 23 + \dl) j + \frac 23} \\
& \les N^{-\frac 12 } \too 0, 
\end{split}
\label{Z6}
\end{align}

\noi
as $N \to \infty$.

Finally, putting 
\eqref{Z2}, \eqref{Z3}, \eqref{Z5}, and \eqref{Z6} together
with $t = N^{-2+\dl}$ and $K = c_0 \log \log N$ for some small $\dl, c_0 > 0$
and proceeding as in \eqref{final}, 
we conclude the proof of Theorem \ref{THM:1}
at the endpoint case $s = -\frac 23$.

\begin{ackno}\rm
 The authors would like to thank the anonymous referee for the helpful comments which improved the quality of the paper.
 I.C.~was supported by the EPSRC New Investigator Award EP/X015688/1.
 T.O.~was supported by the 
 European Research Council (grant no.~864138 ``SingStochDispDyn"). 
Y.W.~was supported by 
 the EPSRC New Investigator Award 
 (grant no.~EP/V003178/1).
T.O.~would like to thank 
the  Centre de recherches math\'ematiques,  Canada, 
for its hospitality, 
where this manuscript was prepared.

\end{ackno}

\end{document}